# On adaptive Bayesian inference

**Yang Xing**

*Centre of Biostochastics*
*Swedish University of Agricultural Sciences*
*SE-901 83 Umeå, Sweden*
*e-mail:* yang.xing@sekon.slu.se

**Abstract:** We study the rate of Bayesian consistency for hierarchical priors consisting of prior weights on a model index set and a prior on a density model for each choice of model index. Ghosal, Lember and Van der Vaart [2] have obtained general in-probability theorems on the rate of convergence of the resulting posterior distributions. We extend their results to almost sure assertions. As an application we study log spline densities with a finite number of models and obtain that the Bayes procedure achieves the optimal minimax rate $n^{-\gamma/(2\gamma+1)}$ of convergence if the true density of the observations belongs to the Hölder space $C^\gamma[0,1]$. This strengthens a result in [1; 2]. We also study consistency of posterior distributions of the model index and give conditions ensuring that the posterior distributions concentrate their masses near the index of the best model.

## 1. Introduction

Selection of models plays a key role in theory of density estimation. Given a collection of models, from the Bayesian point of view it is natural to put a prior on model index and let the resulting posteriors determine a correct model. A rate-adaptive posterior achieves the rate of convergence provided by the best single model from the collection. This paper handles adaptation for density estimation within the Bayesian framework. Suppose that we observe a random sample $X_1, X_2, \ldots, X_n$ generated from a probability distribution $P_0$ with a density function $f_0$ with respect to some dominated $\sigma$-finite measure on a measurable space $\mathbb{X}$. Let $I_n$ denote an at most countable index set for each positive integer $n$. For $\gamma \in I_n$, $\mathcal{P}_{n,\gamma}$ stands for a subset of the density space $\mathbb{F}$ equipped with a $\sigma$-field such that the mapping $(x,f) \mapsto f(x)$ is measurable relative to the product $\sigma$-field on $\mathbb{X} \times \mathcal{P}_{n,\gamma}$. Let $\Pi_{n,\gamma}$ be a probability measure on $\mathcal{P}_{n,\gamma}$ and let $\{\lambda_{n,\gamma} : \gamma \in I_n\}$ be a discrete probability measure on $I_n$. One can therefore define an overall prior $\Pi_n$ with support on $\cup_{\gamma \in I_n} \mathcal{P}_{n,\gamma} \subset \mathbb{F}$ by

$$\Pi_n = \sum_{\gamma \in I_n} \lambda_{n,\gamma} \, \Pi_{n,\gamma}.$$

The corresponding posterior distribution $\Pi_n(\cdot \mid X_1, X_2, \ldots, X_n)$ is a random probability measure with the expression

$$\Pi_n(A \mid X_1, X_2, \ldots, X_n) = \frac{\int_A \prod_{i=1}^n f(X_i)\, \Pi_n(df)}{\int_\mathbb{F} \prod_{i=1}^n f(X_i)\, \Pi_n(df)} = \frac{\int_A R_n(f)\, \Pi_n(df)}{\int_\mathbb{F} R_n(f)\, \Pi_n(df)}$$

$$= \frac{\sum_{\gamma \in I_n} \lambda_{n,\gamma} \int_{\mathcal{P}_{n,\gamma} \cap A} R_n(f)\, \Pi_{n,\gamma}(df)}{\sum_{\gamma \in I_n} \lambda_{n,\gamma} \int_{\mathcal{P}_{n,\gamma}} R_n(f)\, \Pi_{n,\gamma}(df)}$$

for all measurable subsets $A \subset \mathbb{F}$, where $R_n(f) = \prod_{i=1}^n \{f(X_i)/f_0(X_i)\}$ denotes the likelihood ratio. The posterior distribution $\Pi_n(\cdot \mid X_1, X_2, \ldots, X_n)$ is said to be consistent almost surely (or in probability) at a rate at least $\varepsilon_n$ if there exists a constant $r > 0$ such that $\Pi_n(f : d(f, f_0) \geq r\varepsilon_n \mid X_1, X_2, \ldots, X_n) \longrightarrow 0$ almost surely (or in probability) as $n \to \infty$. Throughout this paper we assume that $d$ is a distance bounded above by the Hellinger distance and $d(f, f_1)^s$ is a convex function of $f$ for some positive constant $s$ and any fixed $f_1$. Almost sure convergence and convergence in probability should be understood as to be with respect to the infinite product distribution $P_0^\infty$ of the true distribution $P_0$.

The purpose is to deal with the following problem: assume that for a given density $f_0$ there exists a best model $\mathcal{P}_{n,\beta_n}$ equipped with a prior $\Pi_{n,\beta_n}$ such that the optimal posterior rate is $\varepsilon_{n,\beta_n}$. Find conditions ensuring that the posterior distributions of the hierarchical prior $\Pi_n$ achieve the same rate of convergence as we only use the best single model $\Pi_{n,\beta_n}$ for this $f_0$. Ghosal, Lember and Van der Vaart [1; 2] have studied adaptation to general models and obtained in-probability results on convergence rate. See also Huang [4] and Lember and van der Vaart [5] for related work on Bayesian adaptation. When applying to log spline density models, Theorem 2.1 of [2] leads to adaptation up to a logarithmic factor and it was shown in [2] that the additional logarithmic factor in the convergence rate can be removed by choosing special prior weights $\lambda_{n,\gamma}$ when $I_n$ are finite sets or the priors $\Pi_{n,\gamma}$ are discrete. Our main goal in present paper is to extend work of Ghosal et al. [2] and establish the corresponding almost sure assertions. With an application of our theorems to log spline densities with finitely many models, we successfully take away the logarithmic factor without using the special prior weights $\lambda_{n,\gamma}$ and hence for a true density in $C^\gamma[0,1]$ the posteriors attain the optimal rate of convergence in the minimax sense, which is well known to be $n^{-\gamma/(2\gamma+1)}$. This strengthens Theorem 5.2 in [2] and Theorem 2 in [1]. A related problem is model selection, for which we establish an almost sure result on consistency of posterior distributions of the model index.

We shall use the Hellinger distance $H(f, g) = ||\sqrt{f} - \sqrt{g}||_2$ and its modification $H_*(f, g) = ||(\sqrt{f} - \sqrt{g})(\frac{2}{3}\sqrt{\frac{f}{g}} + \frac{1}{3})^{1/2}||_2$, where $||f||_p = \left(\int_\mathbb{X} |f(x)|^p\, \mu(dx)\right)^{1/p}$. Observe that $H_*(f, g) \neq H_*(g, f)$, see [9] for properties of $H_*(f, g)$. Denote

$$W_{n,\gamma}(\varepsilon) = \{f \in \mathcal{P}_{n,\gamma} : H_*(f_0, f) \leq \varepsilon\},$$

$$A_{n,\gamma}(\varepsilon) = \{f \in \mathcal{P}_{n,\gamma} : d(f_0, f) \leq \varepsilon\}.$$

Throughout this paper the notation $a \lesssim b$ means that $a \leq Cb$ for some positive constant $C$ which is universal or fixed in the proof. Write $a \approx b$ if $a \lesssim b$ and $b \lesssim a$. For a measure $P$ and an integrable function $f$ on $\mathbb{X}$, we let $Pf$ stand for the integral of $f$ on $\mathbb{X}$ with respect to $P$. The notation $N(\delta, \mathcal{G}, d)$ stands for the minimal number of balls of radius $\delta$ relative to the distance $d$ needed to cover a subset $\mathcal{G}$ of $\mathbb{F}$.

## 2. Adaptation and Model Selection

Denote by $\varepsilon_{n,\gamma}$ the usual optimal convergence rate of posteriors by using the single model $\mathcal{P}_{n,\gamma}$ with the prior $\Pi_{n,\gamma}$. We shall use a partition $I_n = I_n^1 \cup I_n^2$ with

$$I_n^1 = \{\gamma \in I_n : \varepsilon_{n,\gamma} \leq \sqrt{H} \varepsilon_{n,\beta_n}\} \quad \text{and} \quad I_n^2 = \{\gamma \in I_n : \varepsilon_{n,\gamma} > \sqrt{H} \varepsilon_{n,\beta_n}\},$$

where $H$ is a fixed constant $\geq 1$.

**Theorem 2.1.** *Suppose that there exist positive constants $H \geq 1$, $E_\gamma$, $\mu_{n,\gamma}$, $G$, $J$, $L$, $C$ and $0 < \alpha < 1$ such that $1 - \alpha > 18\alpha L$, $n\varepsilon_{n,\beta_n}^2 \geq (1 + \frac{1}{C})\log n$, $\sup_{\gamma \in I_n^1} E_\gamma \varepsilon_{n,\gamma}^2 \leq G\varepsilon_{n,\beta_n}^2$, $\sup_{\gamma \in I_n^2} E_\gamma \leq G$ and $\sum_{\gamma \in I_n} \mu_{n,\gamma}^\alpha = O(e^{Jn\varepsilon_{n,\beta_n}^2})$. Let $r$ be a constant with $r \geq \frac{18(C+J+G+3\alpha+2\alpha C)}{1-\alpha-18\alpha L} + \sqrt{H} + 1$ such that*

(1) $N\left(\frac{\varepsilon}{3}, A_{n,\gamma}(2\varepsilon), d\right) \leq e^{E_\gamma n\varepsilon_{n,\gamma}^2}$ *for all* $\gamma \in I_n$ *and* $\varepsilon \geq \varepsilon_{n,\gamma}$,

(2) $\dfrac{\lambda_{n,\gamma} \Pi_{n,\gamma}\left(A_{n,\gamma}(j\varepsilon_{n,\gamma})\right)}{\lambda_{n,\beta_n} \Pi_{n,\beta_n}\left(W_{n,\beta_n}(\varepsilon_{n,\beta_n})\right)} \leq \mu_{n,\gamma} e^{Lj^2 n\varepsilon_{n,\gamma}^2}$ *for all* $\gamma \in I_n^2$ *and* $j \geq r$,

(3) $\dfrac{\lambda_{n,\gamma} \Pi_{n,\gamma}\left(A_{n,\gamma}(j\varepsilon_{n,\beta_n})\right)}{\lambda_{n,\beta_n} \Pi_{n,\beta_n}\left(W_{n,\beta_n}(\varepsilon_{n,\beta_n})\right)} \leq \mu_{n,\gamma} e^{Lj^2 n\varepsilon_{n,\beta_n}^2}$ *for all* $\gamma \in I_n^1$ *and* $j \geq r$,

(4) $\sum_{n=1}^{\infty} \sum_{\gamma \in I_n^2} \dfrac{\lambda_{n,\gamma} \Pi_{n,\gamma}\left(A_{n,\gamma}(r\varepsilon_{n,\gamma})\right) e^{(3+2C)n\varepsilon_{n,\beta_n}^2}}{\lambda_{n,\beta_n} \Pi_{n,\beta_n}\left(W_{n,\beta_n}(\varepsilon_{n,\beta_n})\right)} < \infty.$

*Then*

$$\Pi_n\left(f : d(f, f_0) \geq r\varepsilon_{n,\beta_n} \,\big|\, X_1, X_2, \ldots, X_n\right) \longrightarrow 0$$

*almost surely as* $n \to \infty$.

Clearly, it is enough to assume that all inequalities in Theorem 2.1 hold for all sufficiently large $n$. As a direct consequence of Theorem 2.1, we have

**Corollary 2.2.** *Suppose that there exist positive constants $H \geq 1$, $E_\gamma$, $\mu_{n,\gamma}$, $G$, $J$, $L$, $C$, $F$ and $0 < \alpha < 1$ such that $1 - \alpha > 18\alpha L$, $n\varepsilon_{n,\beta_n}^2 \geq (1 + \frac{1}{C})\log n$, $\sup_{\gamma \in I_n^1} E_\gamma \varepsilon_{n,\gamma}^2 \leq G\varepsilon_{n,\beta_n}^2$, $\sup_{\gamma \in I_n^2} E_\gamma \leq G$ and $\sum_{\gamma \in I_n} \mu_{n,\gamma}^\alpha = \mathrm{O}(e^{Jn\varepsilon_{n,\beta_n}^2})$. Let $r$ be a constant such that $r \geq \frac{18(C+J+G+3\alpha+2\alpha C)}{1-\alpha-18\alpha L} + \sqrt{H} + 1$ and*

(1) $N(\frac{\varepsilon}{3}, A_{n,\gamma}(2\varepsilon), d) \leq e^{E_\gamma n \varepsilon_{n,\gamma}^2}$ *for all $\gamma \in I_n$ and $\varepsilon \geq \varepsilon_{n,\gamma}$,*

(2) $\frac{\lambda_{n,\gamma}}{\lambda_{n,\beta_n}} \leq \mu_{n,\gamma} e^{(L-F)n(\varepsilon_{n,\gamma}^2 \vee \varepsilon_{n,\beta_n}^2)}$ *for all $\gamma \in I_n$,*

(3) $\sum_{\gamma \in I_n^2} \frac{\lambda_{n,\gamma}}{\lambda_{n,\beta_n}} \Pi_{n,\gamma}(A_{n,\gamma}(r\varepsilon_{n,\gamma})) = \mathrm{O}\left(e^{-(3+3C+F)n\varepsilon_{n,\beta_n}^2}\right),$

(4) $\Pi_{n,\beta_n}(W_{n,\beta_n}(\varepsilon_{n,\beta_n})) \geq e^{-Fn\varepsilon_{n,\beta_n}^2}.$

*Then*
$$\Pi_n\big(f : d(f, f_0) \geq r\varepsilon_{n,\beta_n} \,\big|\, X_1, X_2, \ldots, X_n\big) \longrightarrow 0$$
*almost surely as $n \to \infty$.*

Condition (3) of Theorem 2.1 leads adaptation up to a logarithmic factor for log spline density models, see [2] for the corresponding in-probability assertion. The following theorem is useful to remove the logarithmic factor in some cases.

**Theorem 2.3.** *Theorem 2.1 holds for $r \geq \frac{18(C+J+G+3\alpha K+2\alpha CK)}{1-\alpha-18\alpha L} + \sqrt{H} + 1$ if the condition (3) of Theorem 2.1 is replaced by the condition that there exists a constant $K \geq 1$ independent of $n$, $\gamma$, $j$ such that*

(3′) $\frac{\Pi_{n,\gamma}(A_{n,\gamma}(j\varepsilon_{n,\beta_n}))}{\Pi_{n,\gamma}(W_{n,\gamma}(K\varepsilon_{n,\beta_n}))} \leq \mu_{n,\gamma} e^{Lj^2 n \varepsilon_{n,\beta_n}^2}$ *for all $\gamma \in I_n^1$ and $j \geq r$.*

Now we consider the rate of convergence of posterior distributions of the index parameter $\gamma$. Given a subset $I$ of $I_n$, Ghosal et al.[2] introduced the posteriors

$$\Pi_n\big(I \,\big|\, X_1, X_2, \ldots, X_n\big) = \frac{\sum_{\gamma \in I} \lambda_{n,\gamma} \int_{\mathcal{P}_{n,\gamma}} R_n(f) \, \Pi_{n,\gamma}(df)}{\sum_{\gamma \in I_n} \lambda_{n,\gamma} \int_{\mathcal{P}_{n,\gamma}} R_n(f) \, \Pi_{n,\gamma}(df)}.$$

Clearly, the result of Theorem 2.1 implies that

$$\Pi_n\big(\gamma \in I_n : d(f_0, \mathcal{P}_{n,\gamma}) \geq r\varepsilon_{n,\beta_n} \,\big|\, X_1, X_2, \ldots, X_n\big) \longrightarrow 0$$

almost surely as $n \to \infty$. Moreover, we have

**Theorem 2.4.** *Under the same assumptions of Theorem 2.1, we have that*

$$\Pi_n\big(I_n^2 \,\big|\, X_1, X_2, \ldots, X_n\big) \longrightarrow 0$$

almost surely as $n \to \infty$. If furthermore for $I_n^3 = \{\gamma \in I_n : \sqrt{H}\varepsilon_{n,\gamma} < \varepsilon_{n,\beta_n}\}$ we have that

$$\sum_{n=1}^{\infty} \sum_{\gamma \in I_n^3} \frac{\lambda_{n,\gamma} \, \Pi_{n,\gamma}\big(A_{n,\gamma}(r\varepsilon_{n,\beta_n})\big) \, e^{(3+2C)n\varepsilon_{n,\beta_n}^2}}{\lambda_{n,\beta_n} \, \Pi_{n,\beta_n}\big(W_{n,\beta_n}(\varepsilon_{n,\beta_n})\big)} < \infty,$$

then

$$\Pi_n\big(\gamma \in I_n : \tfrac{1}{\sqrt{H}}\varepsilon_{n,\beta_n} \leq \varepsilon_{n,\gamma} \leq \sqrt{H}\varepsilon_{n,\beta_n} \,\big|\, X_1, X_2, \ldots, X_n\big) \longrightarrow 1$$

almost surely as $n \to \infty$.

Since $H$ is an arbitrarily given constant bigger than 1, Theorem 2.4 states that the posterior distributions of model index concentrate their masses on the indices of those models which have approximately the same convergence rate as the correct rate $\varepsilon_{n,\beta_n}$. So Theorem 2.4 can be considered as a general convergence theorem on posterior distributions of model index.

In the situation that there are only two models, one can use the Bayes factor to describe behavior of the posterior of the model index, see [2]. Denote by $BF_n$ the Bayes factor, that is,

$$FB_n := \frac{\lambda_{n,2} \int_{\mathcal{P}_{n,2}} R_n(f) \, \Pi_{n,2}(df)}{\lambda_{n,1} \int_{\mathcal{P}_{n,1}} R_n(f) \, \Pi_{n,1}(df)} = \frac{\Pi_n\big(\{2\} \,\big|\, X_1, X_2, \ldots, X_n\big)}{\Pi_n\big(\{1\} \,\big|\, X_1, X_2, \ldots, X_n\big)}.$$

**Corollary 2.5.** *Suppose that condition (1) of Theorem 2.1 holds and that $\varepsilon_{n,1} > \varepsilon_{n,2} \geq \sqrt{(1+1/C)(\log n)/n}$ for all $n$ and some $C > 0$. Let $r > 700(2C+G+2)$.*
  (i) *If $\Pi_{n,2}\big(W_{n,2}(\varepsilon_{n,2})\big) \geq e^{-n\varepsilon_{n,2}^2}$, $\frac{\lambda_{n,1}}{\lambda_{n,2}} \Pi_{n,1}\big(A_{n,1}(r\varepsilon_{n,1})\big) = \mathrm{O}\big(e^{-(4+3C)n\varepsilon_{n,2}^2}\big)$ and $\frac{\lambda_{n,1}}{\lambda_{n,2}} \leq e^{n\varepsilon_{n,1}^2}$, then $BF_n \to \infty$ almost surely.*
  (ii) *If $\Pi_{n,1}\big(W_{n,1}(\varepsilon_{n,1})\big) \geq e^{-n\varepsilon_{n,1}^2}$, $\frac{\lambda_{n,2}}{\lambda_{n,1}} \Pi_{n,2}\big(A_{n,2}(r\varepsilon_{n,1})\big) = \mathrm{O}\big(e^{-(4+3C)n\varepsilon_{n,1}^2}\big)$ and $\frac{\lambda_{n,2}}{\lambda_{n,1}} \leq e^{n\varepsilon_{n,1}^2}$, then $BF_n \to 0$ almost surely.*

*Proof.* Take $H = J = F = 1$, $L = 2$ and $\alpha = 1/38$. Then $1 - \alpha > 18\alpha L$.

(i) Let $\beta_n = 2$. Then $I_n^1 = \{2\}$ and $I_n^2 = \{1\}$. It follows then from the first assertion of Theorem 2.4 that the denominator of the Bayes factor $BF_n$ tends to zero almost surely as $n \to \infty$ and hence $BF_n \to \infty$ almost surely.

(ii) Let $\beta_n = 1$. Then $I_n^1 = \{1, 2\}$, $I_n^2 = \emptyset$ and $I_n^3 = \{2\}$. It follows then from the second assertion of Theorem 2.4 that the numerator of the Bayes factor $BF_n$ tends to zero almost surely as $n \to \infty$ and hence $BF_n \to 0$ almost surely. The proof of Corollary 2.5 is complete. □

## 3. Log Spline Density Models

Log spline density models were introduced by Stone [7] in his study of sieved maximum likelihood estimators, and were developed by Ghosal, Ghosh and Van der Vaart [3] to Bayesian estimators. Assume that $\big[(k-1)/K_n, k/K_n\big)$ with $k = 1, 2, \ldots, K_n$ is a given partition of the half open interval $[0, 1)$. The space of splines of order $q$ relative to this partition is the set of all functions $f : [0, 1] \mapsto \mathbb{R}$ such that $f$ is $q - 2$ times continuously differentiable on $[0, 1)$ and the restriction of $f$ on each $\big[(k-1)/K_n, k/K_n\big)$ is a polynomial of degree strictly less then $q$. Given $\gamma > 0$, denote $J_{n,\gamma} = q + K_n - 1$ where $q$ is a fixed constant $\geq \gamma$. The space of splines is a $J_{n,\gamma}$-dimensional vector space with a basis $B_1(x), B_2(x), \ldots, B_{J_{n,\gamma}}(x)$ of B-splines, which is a uniformly bounded nonnegative function supported on some interval of length $q/K_n$, see [3] for the details of such a B-spline basis. Assume throughout that the true density $f_0(x) := f_{\theta_0}(x)$ is bounded away from zero and infinity in $[0, 1]$. We consider the $J_{n,\gamma}$-dimensional exponential subfamily of $C^\gamma[0, 1]$ of the form

$$f_\theta(x) = \exp\Big(\sum_{j=1}^{J_{n,\gamma}} \theta_j B_j(x) - c(\theta)\Big),$$

where $\theta = (\theta_1, \theta_2, \ldots, \theta_{J_{n,\gamma}}) \in \Theta_0 = \{(\theta_1, \theta_2, \ldots, \theta_{J_{n,\gamma}}) \in \mathbb{R}^{J_{n,\gamma}} : \sum_{j=1}^{J_{n,\gamma}} \theta_j = 0\}$ and the constant $c(\theta)$ is chosen such that $f_\theta(x)$ is a density in $[0, 1]$. Each prior $\overline{\Pi}_{n,\gamma}$ on $\Theta_0$ induces naturally a prior $\Pi_{n,\gamma}$ on the density set $\mathcal{P}_{n,\gamma} := \{f_\theta(x) : \theta \in \Theta_0\}$. Assume that $J_{n,\gamma} \approx K_n \approx n^{1/(2\gamma+1)}$ and assume that the prior $\overline{\Pi}_{n,\gamma}$ for $\Theta_0$ is supported on $[-M, M]^{J_{n,\gamma}}$ for some $M \geq 1$ and has a density function with respect to the Lebesgue measure, which is bounded on $[-M, M]^{J_{n,\gamma}}$ below by $d^{J_{n,\gamma}}$ and above by $D^{J_{n,\gamma}}$ for two fixed constants $d$ and $D$ with $0 < d \leq D < \infty$. Write $||\theta||_p = \big(\sum_{j=1}^{J_{n,\gamma}} |\theta_j|^p\big)^{1/p}$ and $||f_\theta(x)||_p = \big(\int f_\theta(x)^p dx\big)^{1/p}$ for $1 \leq p \leq \infty$. Take constants $\overline{C}_1 \geq \underline{C}_1 > 0$ such that $\underline{C}_1 ||\theta||_\infty \leq ||\log f_\theta(x)||_\infty \leq \overline{C}_1 ||\theta||_\infty$ for all $\theta \in \Theta_0$, see Lemma 7.3 in [2] for existence of $\underline{C}_1$ and $\overline{C}_1$. Hence $e^{-\overline{C}_1 M} \leq f_\theta(x) \leq e^{\overline{C}_1 M}$ for all $\theta \in \Theta_0$ with $||\theta||_\infty \leq M$. Ghosal et al.([3], Theorem 4.5) proved that, if $f_0 \in C^\gamma[0, 1]$ with $q \geq \gamma \geq 1/2$ and $||\log f_0(x)||_\infty \leq \underline{C}_1 M/2$, the posteriors are consistent in probability at the rate $n^{-\gamma/(2\gamma+1)}$. This result has been strengthened by Xing [9] to the almost sure consistency of the posteriors.

For given priors $\Pi_{n,\gamma}$ on densities and a discrete prior $\{\lambda_{n,\gamma}\}$ on regularity parameters $\gamma$, we get an overall prior $\Pi_n$ on densities as before. Under mild conditions, Ghosal et al. [2] obtained an in-probability theorem on adaptation up to a logarithmic factor for the posteriors. They also showed in [1; 2] that the logarithmic factor can be removed by choosing special prior weights $\lambda_{n,\gamma}$ either when $I_n$ are finite sets or when all the priors $\Pi_{n,\gamma}$ are discrete. Now, for finite index sets $I_n$, we can take away the logarithmic factor without using the special prior weights $\lambda_{n,\gamma}$ and our result moreover is an almost sure statement.

Following [2], we consider prior weights $\lambda_{n,\gamma} = \lambda_\gamma > 0$ for all $n$ and $\gamma \in I_n := \{\gamma \in \mathbb{Q}^+ : \gamma \geq \gamma_1\}$, where $\gamma_1$ is a known positive constant strictly bigger than $1/2$. Now we prove

**Theorem 3.1.** *Let $I_n = \{\gamma_1, \gamma_2, \ldots, \gamma_N\}$ and $\varepsilon_{n,\gamma} = n^{-\gamma/(2\gamma+1)}$ for all $\gamma \in I_n$. If $f_0 \in C^\beta[0, 1]$ with some $\beta \in I_n$ and $||\log f_0(x)||_\infty \leq \underline{C}_1 M$, then for all large*

*constants r,*

$$\Pi_n\{f_\theta : ||f_\theta - f_0||_2 \geq r\varepsilon_{n,\beta} | X_1, \ldots, X_n\} \longrightarrow 0$$

*almost surely as $n \to \infty$.*

*Proof.* We shall apply Theorem 2.3 for the Hellinger distance to the proof. Observe first that $n\varepsilon_{n,\beta}^2 = n^{1/(2\beta+1)} \geq (1 + 1/C) \log n$ when $n$ is large enough and $C = 1$. Take $\mu_{n,\gamma} = \lambda_\gamma/\lambda_\beta$. Conditions (1) of Theorem 2.3 has been verified in [2]. Denote

$$\Theta_{0,M} = \{\theta \in \Theta_0 : ||\theta||_\infty \leq M\},$$
$$C_{J_{n,\gamma}}(\varepsilon) = \{f_\theta : H(f_\theta, f_0) \leq \varepsilon \text{ and } \theta \in \Theta_{0,M}\},$$
$$W_{J_{n,\gamma}}(\varepsilon) = \{f_\theta : H_*(f_\theta, f_0) \leq \varepsilon \text{ and } \theta \in \Theta_{0,M}\}.$$

Since $f_0/f_\theta$ are uniformly bounded above by $e^{(\overline{C}_1+\underline{C}_1)M}$ and below by $e^{-(\overline{C}_1+\underline{C}_1)M}$ for all $\theta \in \Theta_{0,M}$, we have

$$W_{J_{n,\gamma}}(\varepsilon/B) \subset C_{J_{n,\gamma}}(\varepsilon) \subset W_{J_{n,\gamma}}(B\varepsilon)$$

for $B = e^{(\overline{C}_1+\underline{C}_1)M/2}$. Hence, applying Lemma 7.6 and Lemma 7.8 in [2], one can find four positive constants $\underline{A}_1, \underline{A}_2, \overline{A}_1$ and $\overline{A}_2$ such that for all large $n$ and all $\varepsilon > 0$,

$$\Pi_{n,\gamma}(C_{J_{n,\gamma}}(\varepsilon)) \leq \overline{\Pi}_{n,\gamma}(\theta \in \Theta_{0,M} : ||\theta - \theta_{J_{n,\gamma}}||_2 \leq \underline{A}_1\sqrt{J_{n,\gamma}}\,\varepsilon)$$
$$\lesssim (\underline{A}_2\sqrt{J_{n,\gamma}}\,\varepsilon)^{J_{n,\gamma}}$$

and

$$\Pi_{n,\beta}(W_{J_{n,\beta}}(\varepsilon_{n,\beta})) \geq \overline{\Pi}_{n,\beta}(\theta \in \Theta_{0,M} : ||\theta - \theta_{J_{n,\beta}}||_2 \leq \overline{A}_1\sqrt{J_{n,\beta}}\,\varepsilon_{n,\beta})$$
$$\gtrsim (\overline{A}_2\sqrt{J_{n,\beta}}\,\varepsilon_{n,\beta})^{J_{n,\beta}},$$

where $\overline{\Pi}_{n,\gamma}$ is the corresponding prior of $\Pi_{n,\gamma}$ and $\theta_{J_{n,\gamma}}$ minimizes the map $\theta \mapsto H(f_\theta, f_0)$ over $\Theta_{0,M}$. In fact, Lemma 7.6 of [2] yields the first inequality for $0 < \varepsilon < 1/\underline{A}_1$. However, since $||\theta||_\infty \leq M$ for $\theta \in \Theta_{0,M}$ and $J_{n,\gamma} \to \infty$ as $n \to \infty$, the inequality holds even for all $\varepsilon \geq 1/\underline{A}_1$ and large $n$. It then follows from $\sqrt{J_{n,\gamma}}\,\varepsilon_{n,\gamma} \approx n^{1/2(2\gamma+1)}\,n^{-\gamma/(2\gamma+1)} = n^{(1-2\gamma)/2(2\gamma+1)}$ for all $\gamma \in I_n$ and hence for $\gamma = \beta$ that

$$\frac{\Pi_{n,\gamma}(C_{J_{n,\gamma}}(j\varepsilon_{n,\gamma}))}{\Pi_{n,\beta}(W_{J_{n,\beta}}(\varepsilon_{n,\beta}))} \lesssim \frac{\left(\underline{A}_2\,j\,n^{\frac{1-2\gamma}{2(2\gamma+1)}}\right)^{J_{n,\gamma}}}{\left(\overline{A}_2\,n^{\frac{1-2\beta}{2(2\beta+1)}}\right)^{J_{n,\beta}}} =$$

$$exp\Big(\Big(\frac{(1-2\gamma)J_{n,\gamma}}{2(2\gamma+1)} + \frac{(2\beta-1)J_{n,\beta}}{2(2\beta+1)}\Big)\log n + J_{n,\gamma}\log\underline{A}_2 - J_{n,\beta}\log\overline{A}_2 + J_{n,\gamma}\log j\Big).$$

Now, for $\gamma \in I_n^2$ we have that $\varepsilon_{n,\gamma} > \sqrt{H}\varepsilon_{n,\beta} > \varepsilon_{n,\beta}$ which implies $\gamma < \beta$ and $J_{n,\gamma} \gtrsim H J_{n,\beta}$. Therefore, using $\frac{1-2\gamma}{4(2\gamma+1)} \leq \frac{1-2\gamma_1}{4(2\gamma_1+1)} < 0$ for $\gamma \geq \gamma_1 > 1/2$, we get that for large $n$ the exponent in the right hand side of the last equality does not exceed a constant multiple of the following sum

$$J_{n,\gamma} \log n \Big( \frac{1-2\gamma}{4(2\gamma+1)} + \frac{2\beta-1}{(2\beta+1)} \frac{1}{H} \Big) + J_{n,\gamma} \log j$$

$$\leq J_{n,\gamma} \log n \Big( \frac{1-2\gamma}{4(2\gamma+1)} + \frac{1}{H} \Big) + J_{n,\gamma} \log j$$

$$\leq J_{n,\gamma} \log n \frac{1-2\gamma_1}{5(2\gamma_1+1)} + J_{n,\gamma} \log j \leq J_{n,\gamma} \log j \leq L j^2 n \varepsilon_{n,\gamma}^2,$$

where $L$ is any given positive constant, the second inequality holds for a large constant $H$ depending only on $\gamma_1$, and the last inequality follows from $J_{n,\gamma} \approx n\varepsilon_{n,\gamma}^2$ and $j \geq r$ with a large $r$. Hence we have verified condition (2). Similarly, since $n\varepsilon_{n,\gamma}^2/H > n\varepsilon_{n,\beta}^2$ for $\gamma \in I_n^2$, we have that for some $M_1 > 0$ and large $H, M_2 > 0$,

$$\sum_{n=M_2}^{\infty} \sum_{\gamma \in I_n^2} \frac{\lambda_{n,\gamma} \Pi_{n,\gamma}\big(C_{J_{n,\gamma}}(r\varepsilon_{n,\gamma})\big) e^{(3+2)n\varepsilon_{n,\beta}^2}}{\lambda_{n,\beta} \Pi_{n,\beta}\big(W_{J_{n,\beta}}(\varepsilon_{n,\beta})\big)}$$

$$\lesssim \sum_{n=M_2}^{\infty} \sum_{\gamma \in I_n^2} \frac{\lambda_\gamma}{\lambda_\beta} e^{J_{n,\gamma} \log n \big( \frac{1-2\gamma_1}{5(2\gamma_1+1)} + \frac{\log r}{\log n} \big) + \frac{5}{H} n\varepsilon_{n,\gamma}^2}$$

$$\leq \sum_{n=M_2}^{\infty} \sum_{\gamma \in I_n^2} \frac{\lambda_\gamma}{\lambda_\beta} e^{n\varepsilon_{n,\gamma}^2 \log n \, M_1 \frac{1-2\gamma_1}{6(2\gamma_1+1)} + \frac{5}{H} n\varepsilon_{n,\gamma}^2}$$

$$\leq \sum_{n=M_2}^{\infty} \sum_{\gamma \in I_n^2} \frac{\lambda_\gamma}{\lambda_\beta} e^{n\varepsilon_{n,\gamma}^2 \log n \, M_1 \frac{1-2\gamma_1}{7(2\gamma_1+1)}}$$

$$\leq \sum_{n=M_2}^{\infty} \sum_{\gamma \in I_n} \frac{\lambda_\gamma}{\lambda_\beta} e^{-2\log n} = \sum_{n=M_2}^{\infty} \frac{1}{\lambda_\beta} \frac{1}{n^2} < \infty,$$

which yields condition (4) for $C = 1$. Finally, observe that $\varepsilon_{n,\gamma} \leq \sqrt{H}\varepsilon_{n,\beta_n}$ for all $\gamma \in I_n^1$. Since $I_n^1$ contains at most finitely many indices and $\varepsilon_{n,\gamma}$ is the convergence rate of the model $\mathcal{P}_{n,\gamma}$ for $f_0$, there exists a constant $K_1 \geq 1$ such that $\Pi_{n,\gamma}\big(C_{J_{n,\gamma}}(K_1\varepsilon_{n,\beta})\big) > 0$ for all $\gamma \in I_n^1$ and all large $n$. It then follows from $W_{J_{n,\gamma}}(\varepsilon/B) \subset C_{J_{n,\gamma}}(\varepsilon) \subset W_{J_{n,\gamma}}(B\varepsilon)$ and Lemma 7.6 in [2] that for a large $K > K_1$, $\underline{A}_3$ and $\overline{A}_3$,

$$\frac{\Pi_{n,\gamma}\big(C_{J_{n,\gamma}}(j\varepsilon_{n,\beta})\big)}{\Pi_{n,\gamma}\big(W_{J_{n,\gamma}}(K\varepsilon_{n,\beta})\big)} \lesssim \frac{\big(\underline{A}_3 j \sqrt{J_{n,\gamma}} \varepsilon_{n,\beta}\big)^{J_{n,\gamma}}}{\big(\overline{A}_3 K \sqrt{J_{n,\gamma}} \varepsilon_{n,\beta}\big)^{J_{n,\gamma}}} \leq e^{J_{n,\gamma} \log \frac{\underline{A}_3 j}{\overline{A}_3 K}} \leq e^{L j^2 n \varepsilon_{n,\gamma}^2}$$

for all large $j$ and any given $L > 0$ which yields condition (3′), and therefore by Theorem 2.3 we obtain the required convergence with respect to the Hellinger metric, which in our case is stronger than the convergence with respect to the metric $||\cdot||_2$, since densities $f_\theta$ are uniformly bounded for all $\theta \in \Theta_{0,M}$. The proof of Theorem 3.1 is complete. □

For general countable index sets $I_n$, Theorem 2.1 yields adaptation up to a logarithmic factor.

**Theorem 3.2.** *Assume that $\sum_{\gamma \in I_n} \lambda_\gamma^\alpha < \infty$ for some $0 < \alpha < 1$. Let $\varepsilon_{n,\gamma} = n^{-\gamma/(2\gamma+1)}\sqrt{\log n}$ for all $\gamma \in I_n$. If $f_0 \in C^\beta[0,1]$ with some $\beta \in I_n$ and $||\log f_0(x)||_\infty \leq \underline{C}_1 M$, then for all large constants $r$,*

$$\Pi_n\{f_\theta : ||f_\theta - f_0||_2 \geq r\varepsilon_{n,\beta} | X_1, \ldots, X_n\} \longrightarrow 0$$

*almost surely as $n \to \infty$.*

*Proof.* Completely repeating the proof of Theorem 3.1, we obtain the conditions (1), (2) and (4) of Theorem 2.1. To see condition (3), note that $\varepsilon_{n,\gamma} \leq \sqrt{H}\,\varepsilon_{n,\beta}$ for $\gamma \in I_n^1$ and hence $J_{n,\gamma} \log n \lesssim n\varepsilon_{n,\gamma}^2 \leq Hn\varepsilon_{n,\beta}^2$. Since $\sqrt{J_{n,\gamma}}\,\varepsilon_{n,\beta} \approx n^{\frac{1}{2(2\gamma+1)} - \frac{\beta}{2\beta+1}}$ for all $\gamma$, the proof of Theorem 3.1 yields that for some sufficiently large $H$ and all large $n$,

$$\frac{\Pi_{n,\gamma}(C_{J_{n,\gamma}}(j\varepsilon_{n,\beta}))}{\Pi_{n,\beta}(W_{J_{n,\beta}}(\varepsilon_{n,\beta}))} \lesssim \frac{(\underline{A}_2\, j\, \sqrt{J_{n,\gamma}}\,\varepsilon_{n,\beta})^{J_{n,\gamma}}}{(\overline{A}_2\, \sqrt{J_{n,\beta}}\,\varepsilon_{n,\beta})^{J_{n,\beta}}} \leq \left(j\,n^{\frac{1}{2(2\gamma+1)}}\right)^{J_{n,\gamma}} \left(n^{\frac{\beta}{2\beta+1}}\right)^{J_{n,\beta}}$$

$$\leq (j\,n)^{J_{n,\gamma}} n^{J_{n,\beta}} = e^{J_{n,\gamma} \log j + (J_{n,\gamma} + J_{n,\beta}) \log n} \leq e^{Lj^2 n\varepsilon_{n,\beta}^2}$$

for all large $j$ and any given $L > 0$ which yields condition (3), and hence by Theorem 2.1 we conclude the proof of Theorem 3.2. $\square$

## 4. Appendix

Let $L_\mu$ be the space of all nonnegative integrable functions with the norm $||f||_1$. Write $\log 0 = -\infty$ and $0/0 = 0$. We shall adopt the Hausdorff $\alpha$-entropy introduced by Xing and Ranneby in [10].

**Definition 4.1.** *Let $\alpha \geq 0$ and $\mathcal{G} \subset \mathbb{F}$. For $\delta > 0$, the Hausdorff $\alpha$-entropy $J(\delta, \mathcal{G}, \alpha, \Pi, d)$ of the set $\mathcal{G}$ relative to the prior distribution $\Pi$ and the distance $d$ is defined as*

$$J(\delta, \mathcal{G}, \alpha, \Pi, d) = \log \inf \sum_{j=1}^{N} \Pi(B_j)^\alpha,$$

*where the infimum is taken over all coverings $\{B_1, B_2, \ldots, B_N\}$ of $\mathcal{G}$, where $N$ may take the value $\infty$, such that each $B_j$ is contained in some ball $\{f : d(f, f_j) < \delta\}$ of radius $\delta$ and center at $f_j \in L_\mu$.*

Note that it was proved in [10] that for any $0 \leq \alpha \leq 1$ and $\mathcal{G} \subset \mathbb{F}$,

$$e^{J(\delta, \mathcal{G}, \alpha, \Pi, d)} \leq \Pi(\mathcal{G})^\alpha\, N(\delta, \mathcal{G}, d)^{1-\alpha} \leq N(\delta, \mathcal{G}, d).$$

We begin with a lemma which is essentially given in the proof of Theorem 1 of [9].

**Lemma 4.2.** *Let $0 < \alpha \leq 1$, $\mathcal{G} \subset \mathbb{F}$ and $D_{r\varepsilon} = \{f \in \mathcal{G} : d(f, f_0) \geq r\varepsilon\}$ with $r > 2$ and $\varepsilon > 0$. Then we have*

$$E \left( \int_{D_{r\varepsilon}} R_n(f) \, \Pi(df) \right)^\alpha \leq e^{J(\varepsilon, D_{r\varepsilon}, \alpha, \Pi, d) + \frac{\alpha-1}{2}(r-2)^2 n \varepsilon^2}.$$

*Proof of Lemma 4.2.* Since $E \int_{D_{r\varepsilon}} R_n(f) \, \Pi(df) = \Pi(D_{r\varepsilon}) \leq 1$, it suffices to prove the lemma for $0 < \alpha < 1$. Given a constant $\phi > 1$, by the definition of $J(\varepsilon, D_{r\varepsilon}, \alpha, \Pi, d)$ there exist functions $f_1, f_2, \ldots, f_N$ in $L_\mu$ such that $D_{r\varepsilon} \subset \bigcup_{j=1}^N B_j$, where $B_j = D_{r\varepsilon} \cap \{f : d(f_j, f) < \varepsilon\}$ and $\sum_{j=1}^N \Pi(B_j)^\alpha \leq \phi e^{J(\varepsilon, D_{r\varepsilon}, \alpha, \Pi, d)}$. By shrinking $B_j$ if necessary, we may assume that all the sets $B_j$ are disjoint and nonempty. Taking some $g_j \in B_j$ we get that $d(f_j, f_0) \geq d(g_j, f_0) - d(g_j, f_j) \geq (r-1)\varepsilon$. Write

$$\int_{B_j} R_n(f) \, \Pi_n(df) = \Pi_n(B_j) \prod_{k=0}^{n-1} \frac{f_{kB_j}(X_{k+1})}{f_0(X_{k+1})},$$

where $f_{kB_j}(x) = \int_{B_j} f(x) R_k(f) \, \Pi_n(df) / \int_{B_j} R_k(f) \, \Pi_n(df)$ and $R_0(f) = 1$. The function $f_{kB_j}$ was introduced by Walker [8] and can be considered as the predictive density of $f$ with a normalized posterior distribution, restricted on the set $B_j$. So we have

$$E \left( \int_{D_{r\varepsilon}} R_n(f) \, \Pi(df) \right)^\alpha \leq \sum_{j=1}^N \Pi(B_j)^\alpha E \left( \prod_{k=0}^{n-1} \frac{f_{kB_j}(X_{k+1})^\alpha}{f_0(X_{k+1})^\alpha} \right)$$

$$\leq \phi e^{J(\varepsilon, D_{r\varepsilon}, \alpha, \Pi, d)} \max_{1 \leq j \leq N} E \left( \prod_{k=0}^{n-1} \frac{f_{kB_j}(X_{k+1})^\alpha}{f_0(X_{k+1})^\alpha} \right).$$

Since $d(f, f_j)^s$ is a convex function of $f$ and $d(f, f_j) \leq \varepsilon$, Jensen's inequality implies that $d(f_{kB_j}, f_j) \leq \varepsilon$ for all $k$ and hence $d(f_{kB_j}, f_0) \geq d(f_j, f_0) - d(f_j, f_{kB_j}) \geq (r-2)\varepsilon$. Using $d(f_{kB_j}, f_0) \leq H(f_{kB_j}, f_0)$ and following the same lines as the proof of Theorem 1 in [9], one can get that

$$E \left( \int_{D_{r\varepsilon}} R_n(f) \, \Pi(df) \right)^\alpha \leq \phi e^{J(\varepsilon, D_{r\varepsilon}, \alpha, \Pi, d) + \frac{\alpha-1}{2}(r-2)^2 n \varepsilon^2},$$

which by the arbitrariness of $\phi > 1$ concludes the proof of Lemma 4.2. □

*Proof of Theorem 2.1.* Denote $D(\varepsilon) = \{f : d(f, f_0) \geq \varepsilon\}$. Write

$$\int_{D(r\varepsilon_n, \beta_n)} R_n(f) \, \Pi_n(df) = \sum_{\gamma \in I_n} \lambda_{n,\gamma} \int_{\mathcal{P}_{n,\gamma} \cap D(r\varepsilon_n, \beta_n)} R_n(f) \, \Pi_{n,\gamma}(df)$$

$$= \sum_{\gamma \in I_n^1} \lambda_{n,\gamma} \int_{\mathcal{P}_{n,\gamma} \cap D(r\varepsilon_n, \beta_n)} R_n(f) \, \Pi_{n,\gamma}(df)$$

$$+ \sum_{\gamma \in I_n^2} \lambda_{n,\gamma} \int_{\mathcal{P}_{n,\gamma} \cap D(\frac{r}{\sqrt{H}} \varepsilon_{n,\gamma})} R_n(f) \, \Pi_{n,\gamma}(df)$$

$$+ \sum_{\gamma \in I_n^2} \lambda_{n,\gamma} \int_{\mathcal{P}_{n,\gamma} \cap \{f: r\varepsilon_{n,\beta_n} \leq d(f,f_0) < \frac{r}{\sqrt{H}}\varepsilon_{n,\gamma}\}} R_n(f) \, \Pi_{n,\gamma}(df).$$

Since $0 < \alpha < 1$, it follows from the inequalities $x \leq x^\alpha$ for $0 \leq x \leq 1$ and $(x+y)^\alpha \leq x^\alpha + y^\alpha$ for $x, y \geq 0$ that

$$\Pi_n\big(D(r\varepsilon_{n,\beta_n}) \,\big|\, X_1, X_2, \ldots, X_n\big) = \frac{\int_{D(r\varepsilon_{n,\beta_n})} R_n(f) \, \Pi_n(df)}{\int_{\mathbb{F}} R_n(f) \, \Pi_n(df)}$$

$$\leq \bigg( \sum_{\gamma \in I_n^1} \frac{\lambda_{n,\gamma} \int_{\mathcal{P}_{n,\gamma} \cap D(r\varepsilon_{n,\beta_n})} R_n(f) \, \Pi_{n,\gamma}(df)}{\int_{\mathbb{F}} R_n(f) \, \Pi_n(df)} \bigg)^\alpha$$

$$+ \bigg( \sum_{\gamma \in I_n^2} \frac{\lambda_{n,\gamma} \int_{\mathcal{P}_{n,\gamma} \cap D(\frac{r}{\sqrt{H}} \varepsilon_{n,\gamma})} R_n(f) \, \Pi_{n,\gamma}(df)}{\int_{\mathbb{F}} R_n(f) \, \Pi_n(df)} \bigg)^\alpha$$

$$+ \sum_{\gamma \in I_n^2} \frac{\lambda_{n,\gamma} \int_{\mathcal{P}_{n,\gamma} \cap \{f: r\varepsilon_{n,\beta_n} \leq d(f,f_0) < \frac{r}{\sqrt{H}}\varepsilon_{n,\gamma}\}} R_n(f) \, \Pi_{n,\gamma}(df)}{\int_{\mathbb{F}} R_n(f) \, \Pi_n(df)}$$

$$\leq \sum_{\gamma \in I_n^1} \frac{\lambda_{n,\gamma}^\alpha \Big( \int_{\mathcal{P}_{n,\gamma} \cap D(r\varepsilon_{n,\beta_n})} R_n(f) \, \Pi_{n,\gamma}(df) \Big)^\alpha}{\lambda_{n,\beta_n}^\alpha \Big( \int_{\mathcal{P}_{n,\beta_n}} R_n(f) \, \Pi_{n,\beta_n}(df) \Big)^\alpha}$$

$$+ \sum_{\gamma \in I_n^2} \frac{\lambda_{n,\gamma}^\alpha \Big( \int_{\mathcal{P}_{n,\gamma} \cap D(\frac{r}{\sqrt{H}} \varepsilon_{n,\gamma})} R_n(f) \, \Pi_{n,\gamma}(df) \Big)^\alpha}{\lambda_{n,\beta_n}^\alpha \Big( \int_{\mathcal{P}_{n,\beta_n}} R_n(f) \, \Pi_{n,\beta_n}(df) \Big)^\alpha}$$

$$+ \sum_{\gamma \in I_n^2} \frac{\lambda_{n,\gamma} \int_{\mathcal{P}_{n,\gamma} \cap \{f: r\varepsilon_{n,\beta_n} \leq d(f,f_0) < \frac{r}{\sqrt{H}}\varepsilon_{n,\gamma}\}} R_n(f) \, \Pi_{n,\gamma}(df)}{\lambda_{n,\beta_n} \int_{\mathcal{P}_{n,\beta_n}} R_n(f) \, \Pi_{n,\beta_n}(df)}.$$

From $n\varepsilon_{n,\beta_n}^2 \geq (1+\frac{1}{C})\log n$ it turns out that $\sum_{n=1}^\infty e^{-Cn\varepsilon_{n,\beta_n}^2} \leq \sum_{n=1}^\infty 1/n^{1+C} < \infty$. Hence, by Lemma 1 of [9] and the first Borel-Cantelli Lemma, we have that

$$\int_{\mathcal{P}_{n,\beta_n}} R_n(f) \, \Pi_{n,\beta_n}(df) \geq \Pi_{n,\beta_n}\big(W_{n,\beta_n}(\varepsilon_{n,\beta_n})\big) \, e^{-(3+2C)n\varepsilon_{n,\beta_n}^2}$$

almost surely for all large $n$. Thus, we obtain that

$$\Pi_n\big(D(r\varepsilon_{n,\beta_n}) \,\big|\, X_1, X_2, \ldots, X_n\big)$$

$$\leq \sum_{\gamma \in I_n^1} \frac{\lambda_{n,\gamma}^\alpha \, e^{(3+2C)\alpha n\varepsilon_{n,\beta_n}^2} \Big( \int_{\mathcal{P}_{n,\gamma} \cap D(r\varepsilon_{n,\beta_n})} R_n(f) \, \Pi_{n,\gamma}(df) \Big)^\alpha}{\lambda_{n,\beta_n}^\alpha \, \Pi_{n,\beta_n}\big(W_{n,\beta_n}(\varepsilon_{n,\beta_n})\big)^\alpha}$$

$$+ \sum_{\gamma \in I_n^2} \frac{\lambda_{n,\gamma}^\alpha \, e^{(3+2C)\alpha n \varepsilon_{n,\beta_n}^2} \left( \int_{\mathcal{P}_{n,\gamma} \cap D(\frac{r}{\sqrt{H}}\varepsilon_{n,\gamma})} R_n(f)\, \Pi_{n,\gamma}(df) \right)^\alpha}{\lambda_{n,\beta_n}^\alpha \, \Pi_{n,\beta_n}\!\left( W_{n,\beta_n}(\varepsilon_{n,\beta_n}) \right)^\alpha}$$

$$+ \sum_{\gamma \in I_n^2} \frac{\lambda_{n,\gamma} \, e^{(3+2C)n \varepsilon_{n,\beta_n}^2} \int_{\mathcal{P}_{n,\gamma} \cap \{f:\, r\varepsilon_{n,\beta_n} \le d(f,f_0) < \frac{r}{\sqrt{H}}\varepsilon_{n,\gamma}\}} R_n(f)\, \Pi_{n,\gamma}(df)}{\lambda_{n,\beta_n} \, \Pi_{n,\beta_n}\!\left( W_{n,\beta_n}(\varepsilon_{n,\beta_n}) \right)}$$

$$:= a_n + b_n + c_n$$

almost surely for all large $n$. Given $\delta > 0$, we have

$$P_0^\infty \big\{ \Pi_n\big( D(r\varepsilon_{n,\beta_n}) \,\big|\, X_1, X_2, \ldots, X_n \big) \ge \delta \big\} \le P_0^\infty \{ a_n + b_n + c_n \ge \delta \}$$

$$\le P_0^\infty\{a_n \ge \delta/3\} + P_0^\infty\{b_n \ge \delta/3\} + P_0^\infty\{c_n \ge \delta/3\} \le \frac{3}{\delta} E a_n + \frac{3}{\delta} E b_n + \frac{3}{\delta} E c_n.$$

It turns out from Fubini's theorem and condition (4) that

$$\sum_{n=1}^\infty E c_n$$

$$= \sum_{n=1}^\infty \sum_{\gamma \in I_n^2} \frac{\lambda_{n,\gamma} \, e^{(3+2C)n\varepsilon_{n,\beta_n}^2} \, \Pi_{n,\gamma}\big( \mathcal{P}_{n,\gamma} \cap \{ f:\, r\varepsilon_{n,\beta_n} \le d(f,f_0) < \frac{r}{\sqrt{H}}\varepsilon_{n,\gamma} \} \big)}{\lambda_{n,\beta_n} \, \Pi_{n,\beta_n}\!\left( W_{n,\beta_n}(\varepsilon_{n,\beta_n}) \right)}$$

$$\le \sum_{n=1}^\infty \sum_{\gamma \in I_n^2} \frac{\lambda_{n,\gamma} \, e^{(3+2C)n\varepsilon_{n,\beta_n}^2} \, \Pi_{n,\gamma}\big( A_{n,\gamma}(r\varepsilon_{n,\gamma}) \big)}{\lambda_{n,\beta_n} \, \Pi_{n,\beta_n}\!\left( W_{n,\beta_n}(\varepsilon_{n,\beta_n}) \right)} < \infty.$$

On the other hand, let $[r]$ be the largest integer less than or equal to $r$ and let $D_{n,\gamma,j} = \{ f \in \mathcal{P}_{n,\gamma} :\, j\varepsilon_{n,\beta_n} \le d(f,f_0) < 2j\varepsilon_{n,\beta_n} \}$. Then for any $\gamma \in I_n^1$ we have

$$\mathcal{P}_{n,\gamma} \cap D(r\varepsilon_{n,\beta_n}) \subset \mathcal{P}_{n,\gamma} \cap D([r]\varepsilon_{n,\beta_n}) = \bigcup_{j=[r]}^\infty D_{n,\gamma,j}$$

and hence

$$E a_n \le \sum_{\gamma \in I_n^1} \sum_{j=[r]}^\infty \frac{\lambda_{n,\gamma}^\alpha \, e^{(3+2C)\alpha n \varepsilon_{n,\beta_n}^2} \, E\!\left( \int_{D_{n,\gamma,j}} R_n(f)\, \Pi_{n,\gamma}(df) \right)^\alpha}{\lambda_{n,\beta_n}^\alpha \, \Pi_{n,\beta_n}\!\left( W_{n,\beta_n}(\varepsilon_{n,\beta_n}) \right)^\alpha}.$$

Since $r \ge \sqrt{H}+1$, we have that $j\varepsilon_{n,\beta_n} \ge [r]\varepsilon_{n,\gamma}/\sqrt{H} \ge \varepsilon_{n,\gamma}$ for all $\gamma \in I_n^1$ and $j \ge [r]$. It then follows from Lemma 4.2, Lemma 1 of [10] and condition (1) that

$$E \left( \int_{D_{n,\gamma,j}} R_n(f)\, \Pi_{n,\gamma}(df) \right)^\alpha \le e^{J\left( \frac{j\varepsilon_{n,\beta_n}}{3},\, A_{n,\gamma}(2j\varepsilon_{n,\beta_n}),\, \alpha,\, \Pi_{n,\gamma},\, d \right) + \frac{\alpha-1}{18} j^2 n \varepsilon_{n,\beta_n}^2}$$

$$\le \Pi_{n,\gamma}\big( A_{n,\gamma}(2j\varepsilon_{n,\beta_n}) \big)^\alpha \, N\!\Big( \frac{\varepsilon}{3},\, A_{n,\gamma}(2\varepsilon),\, d \Big)^{1-\alpha} \, e^{\frac{\alpha-1}{18} j^2 n \varepsilon_{n,\beta_n}^2}$$

$$\leq \Pi_{n,\gamma}\big(A_{n,\gamma}(2j\varepsilon_{n,\beta_n})\big)^\alpha e^{E_\gamma n\varepsilon_{n,\gamma}^2 + \frac{\alpha-1}{18}j^2 n\varepsilon_{n,\beta_n}^2}.$$

Thus, by the assumption that $E_\gamma \varepsilon_{n,\gamma}^2 \leq G\varepsilon_{n,\beta_n}^2$ for $\gamma \in I_n^1$, we have

$$Ea_n \leq \sum_{\gamma \in I_n^1} \sum_{j=[r]}^\infty \frac{\lambda_{n,\gamma}^\alpha \, \Pi_{n,\gamma}\big(A_{n,\gamma}(2j\varepsilon_{n,\beta_n})\big)^\alpha e^{(3\alpha+2\alpha C+G+\frac{\alpha-1}{18}j^2)n\varepsilon_{n,\beta_n}^2}}{\lambda_{n,\beta_n}^\alpha \, \Pi_{n,\beta_n}\big(W_{n,\beta_n}(\varepsilon_{n,\beta_n})\big)^\alpha},$$

which by condition (3) does not exceed

$$\sum_{\gamma \in I_n^1} \sum_{j=[r]}^\infty \mu_{n,\gamma}^\alpha \, e^{(3\alpha+2\alpha C+G+\alpha L j^2+\frac{\alpha-1}{18}j^2)n\varepsilon_{n,\beta_n}^2}$$

$$= \mathrm{O}\Big(\sum_{j=[r]}^\infty e^{(J+G+3\alpha+2\alpha C+\alpha L j^2+\frac{\alpha-1}{18}j^2)n\varepsilon_{n,\beta_n}^2}\Big)$$

$$= \mathrm{O}\Big(e^{(J+G+3\alpha+2\alpha C)n\varepsilon_{n,\beta_n}^2} \sum_{j=[r]}^\infty e^{(\alpha L+\frac{\alpha-1}{18})j n\varepsilon_{n,\beta_n}^2}\Big)$$

$$= \mathrm{O}\Big(e^{(J+G+3\alpha+2\alpha C)n\varepsilon_{n,\beta_n}^2} \frac{e^{(\alpha L+\frac{\alpha-1}{18})[r]n\varepsilon_{n,\beta_n}^2}}{1-e^{(\alpha L+\frac{\alpha-1}{18})n\varepsilon_{n,\beta_n}^2}}\Big)$$

$$= \mathrm{O}\Big(e^{(J+G+3\alpha+2\alpha C+\alpha L[r]+\frac{\alpha-1}{18}[r])n\varepsilon_{n,\beta_n}^2}\Big) = \mathrm{O}\Big(e^{-Cn\varepsilon_{n,\beta_n}^2}\Big) = \mathrm{O}\Big(\frac{1}{n^{1+C}}\Big),$$

where the first equality follows from $\sum_{\gamma \in I_n} \mu_{n,\gamma}^\alpha = \mathrm{O}(e^{Jn\varepsilon_{n,\beta_n}^2})$, the third one from $1-\alpha > 18\alpha L$, the next last one from $r \geq \frac{18(C+J+G+3\alpha+2\alpha C)}{1-\alpha-18\alpha L}+1$ and the last one from $n\varepsilon_{n,\beta_n}^2 \geq (1+\frac{1}{C})\log n$. Therefore, we have that $\sum_{n=1}^\infty Ea_n < \infty$. On the other hand, observe that $\varepsilon_{n,\gamma} > \sqrt{H}\varepsilon_{n,\beta_n} \geq \varepsilon_{n,\beta_n}$ for $\gamma \in I_n^2$. So, using the same argument as the above, one can get that

$$Eb_n \leq \sum_{\gamma \in I_n^2} \sum_{j=[r]}^\infty \mu_{n,\gamma}^\alpha \, e^{(3\alpha+2C\alpha+G)n\varepsilon_{n,\beta_n}^2 + (\alpha L j^2+\frac{\alpha-1}{18}j^2)n\varepsilon_{n,\gamma}^2}$$

$$= \mathrm{O}\Big(\sum_{j=[r]}^\infty e^{(J+G+3\alpha+2C\alpha+\alpha L j^2+\frac{\alpha-1}{18}j^2)n\varepsilon_{n,\beta_n}^2}\Big) = \mathrm{O}\Big(\frac{1}{n^{1+C}}\Big),$$

which yields that $\sum_{n=1}^\infty Eb_n < \infty$. Thus, we have proved that

$$\sum_{n=1}^\infty P_0^\infty \big\{\Pi_n\big(f: d(f,f_0) \geq r\varepsilon_{n,\beta_n} \,\big|\, X_1, X_2, \ldots, X_n\big) \geq \delta\big\} < \infty,$$

and then by the first Borel-Cantelli Lemma we get that

$$\Pi_n\big(f: d(f,f_0) \geq r\varepsilon_{n,\beta_n} \,\big|\, X_1, X_2, \ldots, X_n\big) < \delta$$

almost surely for all large $n$. The proof of Theorem 2.1 is complete. $\square$

*Proof of Theorem 2.3.* The proof of Theorem 2.3 is in fact a slight modification of the proof of Theorem 2.1. We only need to repeat the proof of Theorem 2.1 except that we shall apply the following inequalities

$$\bigg(\sum_{\gamma \in I_n^1} \frac{\lambda_{n,\gamma} \int_{\mathcal{P}_{n,\gamma} \cap D(r\varepsilon_{n,\beta_n})} R_n(f)\, \Pi_{n,\gamma}(df)}{\int_{\mathbb{F}} R_n(f)\, \Pi_n(df)}\bigg)^\alpha$$

$$\leq \sum_{\gamma \in I_n^1} \bigg(\frac{\int_{\mathcal{P}_{n,\gamma} \cap D(r\varepsilon_{n,\beta_n})} R_n(f)\, \Pi_{n,\gamma}(df)}{\int_{\mathcal{P}_{n,\gamma}} R_n(f)\, \Pi_{n,\gamma}(df)}\bigg)^\alpha$$

and

$$\int_{\mathcal{P}_{n,\gamma}} R_n(f)\, \Pi_{n,\gamma}(df) \geq \Pi_{n,\gamma}\big(W_{n,\gamma}(K\varepsilon_{n,\beta_n})\big)\, e^{-(3+2C)Kn\varepsilon_{n,\beta_n}^2}.$$

The details of the proof of Theorem 2.3 are therefore omitted. □

*Proof of Theorem 2.4.* The first assertion of Theorem 2.4 follows from the proof of Theorem 2.1. The second assertion follows similarly by applying the partition

$$\mathcal{P}_{n,\gamma} = \bigcup_{\gamma \in I_n^1 \setminus I_n^3} \mathcal{P}_{n,\gamma} \cup \bigcup_{\gamma \in I_n^3} \{f \in \mathcal{P}_{n,\gamma} : d(f,f_0) < r\varepsilon_{n,\beta_n}\}$$

$$\cup \bigcup_{\gamma \in I_n^3} \{f \in \mathcal{P}_{n,\gamma} : d(f,f_0) \geq r\varepsilon_{n,\beta_n}\} \cup \bigcup_{\gamma \in I_n^2} \{f \in \mathcal{P}_{n,\gamma} : d(f,f_0) < \frac{r}{\sqrt{H}}\varepsilon_{n,\gamma}\}$$

$$\cup \bigcup_{\gamma \in I_n^2} \{f \in \mathcal{P}_{n,\gamma} : d(f,f_0) \geq \frac{r}{\sqrt{H}}\varepsilon_{n,\gamma}\}.$$

So we omit the details of the proof of Theorem 2.4. □

## Acknowledgements

It is a great pleasure for me to thank Bo Ranneby for many helpful discussions. I would also like to thank the editor and the associate editor for their helpful comments.